\documentclass[12pt]{article}
\usepackage[dvips]{epsfig}
\usepackage{amsmath}
\usepackage{amsfonts}
\usepackage{amssymb}
\usepackage{graphicx}

\newtheorem{cor}{Corollary}[section]

\numberwithin{equation}{section}

\newcommand{\eps}{\varepsilon}

\newcommand{\M}{\mathcal{M}}

\newcommand{\C}{\mathcal{C}}
\newcommand{\D}{\mathcal{D}}
\renewcommand{\phi}{\varphi}
\newcommand{\F}{\mathcal{F}}

\newcommand{\U}{\mathcal{U}}
\renewcommand{\H}{\mathcal{H}}
\renewcommand{\chi}{\mathcal{X}}

\newcommand{\p}{\partial}
\newcommand{\R}{{\mathbb R}}
\newcommand{\B}{{\mathbb B}}
\newcommand{\noo}{\noindent}
\newcommand{\wtt}{\widetilde}
\newcommand{\lan}{\langle}
\newcommand{\ran}{\rangle}
\newcommand{\Om}{\Omega}
\newcommand{\pom}{\partial\Omega}
\newcommand{\lam}{\lambda}

\begin{document}

\title{\Large \textbf{\boldmath Singularity Profile\\
  in the Mean Curvature Flow }}

\author{{\large Weimin Sheng\ \ \ \ \ \  Xu-Jia Wang }}
\date{\vspace{-8mm}}

\maketitle

\begin{abstract}
In this paper we study the geometry of first time singularities of
the mean curvature flow. By the curvature pinching estimate of
Huisken and Sinestrari, we prove that a mean curvature flow of
hypersurfaces in the Euclidean space $\R^{n+1}$ with positive mean
curvature is $\kappa$-noncollapsing, and a blow-up sequence
converges locally smoothly along a subsequence to a smooth, convex
blow-up solution. As a consequence we obtain a local Harnack
inequality for the mean convex flow.
\end{abstract}

\thispagestyle{empty}


\baselineskip=17.5pt
\parskip=3pt

\section{Introduction}\hskip600pt
\footnote{The first author was supported by NSFC10771189  and
10831008; the second author was supported by ARC DP0664517 and
DP0879422.}
\vskip-20pt

In this paper we are concerned with the geometry of the first time
singularities of the mean curvature flow of closed, smooth
hypersurfaces in the Euclidean space $\R^{n+1}$. The singularity
profile has been studied by many authors [5,8-13,15,21,22], but in
general it can be extremely complicated. For example an open problem
is whether for any $k>1$, there is a self-similar toric solution of
genus $k$ in $\R^3$. Therefore as in [8-13, 21,22], we restrict to
the mean curvature flow of hypersurfaces with positive mean
curvature, namely the {\it mean convex flow}.

For the mean convex flow, Huisken and Sinestrari [11,12] proved the
following crucial one-side curvature pinching estimate, namely at
any point up to the first time singularity,
 \begin{equation}
 \lam_1\ge -\phi(\lam_1+\cdots+\lam_n),
 \end{equation}
where $\lam_1\le \cdots\le\lam_n$ are the principal curvatures and
$\phi$ is a nonnegative function depending on the initial
hypersurface and satisfying $\phi(t)/t\to 0$ as $t\to \infty$. This
estimate corresponds to a one-side curvature pinching estimate of
Hamilton and Ivey (see [6]) for the Ricci flow of 3-manifolds. By
estimate (1.1) and Hamilton's Harnack inequality [7], Huisken and
Sinestrari proved that the blow-up sequence at the maximal mean
curvature points of type II singularities converges along a
subsequence to a convex translating solution.

In this paper we study the asymptotic behavior of general blow-up
sequences. First we prove (see \S3 for terminology and notation).

\vskip5pt

\noo{\bf Theorem 1.1}. {\it Let $\F=\{F_t\}_{t\in [0, T)}$ be a
smooth, mean convex flow of closed hypersurface in $\R^{n+1}$, $n\ge
2$, with first singular time $T$. Then any blow-up sequence
$\F^k=\{F^k_t\}$ at the first time singularity sub-converges locally
smoothly to a convex solution $\F'=\{F'_t\}_{t\in (-\infty, T')}$,
where $0\le T'\le +\infty$.}

\vskip5pt

The local smooth convergence means for any $R>0$, $\F^k\cap Q_R$
converges in the $C^3$ topology, where $Q_R$ denotes the parabolic
cube in spacetime.

We then prove the $\kappa$-noncollapsing. By our definition in \S3,
the $\kappa$-noncollapsing means a pinching of the solution by a
sphere from inside. Recall that (1.1) then means a pinching from
outside.

\vskip5pt

\noo{\bf Theorem 1.2}. {\it Let $\F=\{F_t\}_{t\in [0, T)}$ be a
smooth, mean convex flow of closed hypersurface in $\R^{n+1}$, $n\ge
2$. Then $\F$ is $\kappa$-noncollapsing up to the first singular
time $T$, for some constant $\kappa>0$ depending only on $n$ and the
initial surface $F_0$.}

Note that the $\kappa$-noncollapsing is invariant under dilation.
Hence Theorem 1.2 implies a blow-up sequence at any points in the
flow $\F$ is $\kappa$-noncollapsing with the same $\kappa$. By the
convexity of blow-up solutions, we can improve Theorem 1.2 to

\vskip5pt

\noo{\bf Theorem 1.3}. {\it There exists a constant $\kappa^*>0$
depending only on $n$, such that any normalized limit flow to any
given mean convex flow in $\R^{n+1}$ at the first time singularities
is $\kappa^*$-noncollapsing.}

\vskip5pt

The normalized condition in Theorem 1.3 is only to exclude
hyperplane as limit flows. Theorem 1.3 can be applied to mean
curvature flow with surgeries. To prove Theorem 1.2, we show that a
pair of hyperplanes, including a multiplicity two plane, cannot be a
blow-up solution. From Theorem 1.2, we have

\begin{cor} 
The grim reaper is not blow-up solution.
\end{cor}

Here the grim reaper is the translating solution determined by
hypersurface $\{x\in\R^{n+1}:\ x_{n+1}=-\log \cos x_1\}$, which is
the product of $\R^{n-1}$ with the curve $x_{n+1}=-\log \cos x_1$ in
a 2-plane.

\begin{cor} 
The set of normalized blow-up solutions is compact.
\end{cor}

By our definition in \S3, the mean curvature of a normalized blow-up
solution is equal to 1 at the origin in space time. Hence by Theorem
1.1, a hyperplane is not a normalized blow-up solution. If one
includes hyperplanes as blow-up solutions (at regular points), then
the set of all blow-up solutions whose mean curvature is uniformly
bounded by a positive constant is compact.

From Theorems 1.1 and 1.2, we have a local Harnack inequality.

\begin{cor} 
 Let $\F =\{F_t\}_{t\in [0, T)}$ be a mean convex flow. Then for any $R>0$, there exists
$\delta>0$, depending only on $n, R$, and the initial condition
$F_0$, such that for any point $p=(x_0, t_0)\in \F$ at which the
mean curvature is greater than $\delta^{-1}$, and any point $q\in\F$
in the parabolic cube $Q_{R/H}(p)$, we have
 \begin{equation}
 \delta H(q)\le H(p)\le \delta^{-1} H(q), 
\end{equation}
where $H$ is the mean curvature at $p$.
\end{cor}

Note that if one can prove the local Harnack inequality first, by
the curvature pinching estimate (1.1), one obtains immediately
Theorems 1.1 and 1.2.

Theorem 1.1 asserts that every blow-up solution is convex. But to
understand the geometry of the singularities more precisely, we wish
to classify all blow-up solutions, or more generally all ancient
convex solutions, as a blow-up solution must be ancient. The study
of ancient convex solutions was carried out by the second author in
[20]. The following results have been obtained:
\newline
$\text{\ \ \ }$ (a) A blow-up translating solution in $\R^3$ is
rotationally symmetric.
\newline
$\text{\ \ \ }$ (b) There exists non-rotationally symmetric
translating solution in $\R^{n+1}$ for $n>2$.
\newline
$\text{\ \ \ }$ (c) The blow-down of a blow-up solution is a
shrinking sphere or cylinder.
\newline
$\text{\ \ \ }$ (d) Let $\F=\{F_t\}_{t\in [0, T)}$ be a smooth, mean
convex flow of closed surface in $\R^3$. Then at any point $x_0\in
F_{t_0}$ with large mean curvature, the surface $F_{t_0}$ is, after
normalization, very close to either a cylinder or a convex cap.

Part (d) can be restated as follows: For any $\eps>0$, $R>1$, there
exists $r_{\eps, R}>1$ such that for any normalized blow-up solution
$\F'=\{F'_t\}$ in $\R^3$, there is a compact set $G\subset F'_{0}$
with diameter $\le r_{\eps, R}$, such that for any $p\in
F'_0\backslash G$,  the component $\B_{R/H}(p)\cap F'_0$ is, after
normalization, in the $\eps$-neighborhood of the cylinder
$S^1\times\R^1$, where $H$ is the mean curvature of $F'_0$ at $p$.

\vskip5pt

The above results show that the singularity profile of mean convex
flow is in line with those of Perelman for the Ricci flow of
3-manifolds [17]. In particular our results for mean convex flow in
$\R^3$ correspond to those for the Ricci flow, therefore it is
conceivable that one can carry out surgeries as in [13]. We will not
get into details in this direction. We point out that the idea of
our proof is different from that of Perelman. Perelman needs to
prove the $\kappa$-noncollapsing and classify ancient
$\kappa$-noncollapsing solutions before he obtains the local smooth
convergence of blow-up sequences. We use the curvature pinching
estimate (1.1) and some basic estimates of parabolic equations to
prove the smooth convergence (Theorem 1.1) first, and then use it to
obtain other results.

Some related results have been obtained in a recent paper [13], in
which Huisken and Sinestrari studied the mean curvature flow with
surgeries in $\R^{n+1}$. They proved that when $n\ge 3$ and the sum
of the least two principal curvatures $\lam_1$ and $\lam_2$ is
positive, a blow-up sequence sub-converges locally smoothly to a
convex blow-up solution and every blow-up solution is very close to
a canonical one.

In an earlier work [21] (\S12), prior to Perelman's paper [17],
White proved that the grim reaper and a multiplicity two plane
cannot be blow-up solutions, which was used in [22] to prove that a
blow-up sequence sub-converges locally smoothly to a convex blow-up
solution. Therefore Theorems 1.1 and 1.2 belong to White [21, 22].
His proof contains many novel ideas but nevertheless is very
involved. It uses various results in geometric measure theory and
minimal surface theory. It is difficult for a reader without
background in these areas to understand.

Most recently, Ben Andrews discovered a direct geometric proof of
the $\kappa$-noncollapsing, which, together with (1.1), also implies
the smooth convergence to a convex solution in Theorem 1.1. There
are many other interesting geometric flows in Euclidean space or
Riemannian manifolds. Different methods may apply to different
situations. Our argument may apply to curvature flows of
hypersurfaces provided one can establish a similar one-side
curvature pinching estimate as (1.1) and related interior a priori
estimates.

This paper is arranged as follows. We include in \S2 some
preliminaries. In \S3 we introduce the terminology and notation. In
\S4 we prove the convexity of limit flows. In \S5 we prove that a
semi-noncollapsing blow-up sequence converges locally smoothly to a
limit flow.  In \S6 we prove any blow-up sequence is
semi-noncollapsing. The proof of Theorems 1.1-1.3 is finished in
\S7. Finally in \S8 we discuss tangent flows.

The authors would like to thank Shizhong Du and Oliver Schn\"urer
for carefully reading this paper and their helpful comments.

\newpage
\vskip0pt

\baselineskip=17.0pt

\section{Preliminaries}

A time-dependent, smooth embedding $F_{t}=F(\cdot, t)$:  $\M\to
\R^{n+1}$, where $t\in [0, T)$, is a solution to the mean curvature
flow if
\begin{equation}
\frac{\p}{\p t} F(p, t) =H \cdot\nu  ,
 \ \ p\in\M, \ t\in [0, T),
\end{equation}
where $H$ is the mean curvature, $\nu$ is the unit inward normal.
When $F_t(\M)$ is locally a graph, we have the parabolic equation
\begin{equation}
\frac{\p u}{\p t}=\sqrt{1+\left| Du\right|^{2}}\,
 \text{div}(\frac{Du}{\sqrt{1+\left| Du\right|^{2}}}) .
\end{equation}

It is well known that the mean curvature $H$ satisfies the equation
[8]
\begin{equation}  
 \frac{\p H}{\p t}
  =\Delta H+H | A |^{2} .
 \end{equation}
Since $F_0(\M)$ is a closed hypersurface, it follows that if $H\ge
0$ at $t=0$, then $H>0$ at $t>0$. From (2.3) we also see that if $u$
is a solution of (2.2) in $Q^n_r:=\B^n_r(0)\times [0, r^2]$ with
$0\le H\le C$, then we have the Harnack inequality $H(0, r^2)\ge
C_1\sup H$, where the sup is taken in the domain
$\B^n_{r/2}(0)\times [\frac 14r^2, \frac 34r^2]$. In particular if
$H(0, r^2)=0$, then $H\equiv 0$ in $Q^n_r$, and $u$ is a linear
function by (2.4) below.

For a mean convex flow $\F$, we also have the estimate for the
second fundamental form $A$,
 \begin{equation}
 |A|^2\le CH^2
 \end{equation}
for some constant $C$ depending on $n$ and $F_0$. Indeed, (2.4)
follows from the following equation [10],
$$\frac{\p }{\p t}(\frac{\left| A\right|^{2}}{H^{2}} )
 =\Delta (\frac{\left| A\right|^{2}}{H^{2}})
 + \frac{2}{H}\lan \nabla H,
           \nabla (\frac{\left| A\right|^{2}}{H^{2}})\ran
 -\frac{2}{H^{4}}\left| H\nabla_{i}h_{jk}-\nabla_{i}H\cdot
 h_{jk}\right|^{2}.$$

For the mean curvature flow of graphs, we have the interior gradient
and second derivative estimates. That is if $u$ is a smooth solution
to (2.2) in $Q_r^n$, then for $0<t\le r^2$,
\begin{align}
 &|Du|(0, t)\le \text{exp}(C_1+C_2M^2(\frac {1}{r^2}+\frac 1t)),\\
 &|D^2u|(0, t)\le C_3(1
   +{\sup}_{Q_r^n}|Du|^4(\frac {1}{r^2}+\frac 1t)),
\end{align}
where $M=\sup_{Q_r^n} |u|$, $C_1, C_2$ depend only on $n$, and $C_3$
depends on $n$ and $M$.

Estimates (2.5) and (2.6) were established in [5]. The estimate
(2.5) in [5] also depends on the Lipschitz continuity of the initial
condition. But by choosing the auxiliary function
$$G(x', t, \xi)=t\rho(x')\phi(u)\log u_\xi(x',t)$$
as in [19],  where $x'=(x_1, \cdots, x_n)$ (see notation in \S3),
$\rho(x')=1-|x'|^2/r^2$ is a cut-off function, $\phi(u)=1+u/M$, one
obtains (2.5) by the computation in [19].


\vskip0pt

\section{Terminology and notation}

First we recall the terminology in [21]. Let $\F=\{F_t\}_{t\in [0,
T)}$ be a mean convex flow, which develops first time singularity at
time $T$. For any sequences $p_k=(x_k, t_k)\in\F$ and
$a_k\to\infty$, let $\F^k=\D_{a_k, p_k}(\F)$, where
\begin{equation}
 \D_{a_k, p_k}:\ (x, t)\to \big(a_k(x-x_k), a_k^2(t-t_k)\big) ,\tag {3.1}
 \end{equation}
namely one first makes a translation such that $p_k=(x_k, t_k)$
becomes the origin in space-time, and then make a parabolic dilation
of scale $a_k$. We get a {\it blow-up sequence} $\F^k$.   If
$a_k=H(x_k, t_k)$ is the mean curvature of $\F$ at $(x_k, t_k)$, we
call $\F^k$ a {\it normalized} blow-up sequence at $(x_k, t_k)$.
Note that for a normalized blow-up sequence, $(x_k, t_k)$ must
converges to a singular point (due to the assumption $a_k\to
\infty$).

If $\F^k$ converges locally smoothly (i.e. smoothly in any compact
set) to a mean curvature flow $\F'$, we say $\F'$ is a {\it limit
flow}, or a {\it blow-up solution}. If $\F^k$ is a normalized
blow-up sequence, the limit flow (blow-up solution) is accordingly
called normalized.

{\it Tangent flow}. A limit flow is called tangent flow if the
sequence $(x_k, t_k)=(x_0, t_0)$ is a fixed point.

{\it Ancient solution}. A solution to the mean curvature flow is
ancient if it exists from time $t=-\infty$.

{\it Eternal solution}. A solution to the mean curvature flow is
eternal if it exists for time $t$ from $-\infty$ to $\infty$.

{\it Collapsing}. A blow-up sequence $\F^k$ is collapsing if there
is no ball $\B_r(x_0)$ contained in $U^k_{t}$ (for some $x_0$ and
any $t<0$) for infinitely many $k$'s (see notation below).

{\it Semi-noncollapsing}. A blow-up sequence $\F^k$ is
semi-noncollapsing if there is ball $\B_r(x_0)$ and a time $t_0$
such that $\B_r(x_0)\subset U^k_{t_0}$ for infinitely many $k$'s.

{\it $\kappa$-noncollapsing}. For any point $(x, t)\in\F$, we say
$\F$ is $\kappa$-noncollapsing at $(x, t)\in\F$ if $r_{x, t}H(x,
t)\ge \kappa$, where $\kappa>0$ is a constant, $H(x, t)$ is the mean
curvature of $\F$ at $(x, t)$, and
 \begin{equation}
 r_{x, t}=
 \sup\{\rho:\ \B_\rho (z)\subset U_t\ \text{and}\  |z-x|=\rho\}.
   \tag {3.2}\label {3.2}
   \end{equation}
($r_{x, t}$ is the radius of the largest ball which is contained in
$U_t$ and tangential to $F_t$ at $x$). We say $\F$ is
$\kappa$-noncollapsing if it is $\kappa$-noncollapsing at any point
before the first time singularity.

{\it The grim reaper}:  It is the translating solution given by
$x_{n+1}=-\log \cos x_1$ (it is the product of a curve with
$\R^{n-1}$).

\noo{\bf Remark 3.1.}
\newline
(i) By (3.1), the point $(0, 0)$ always belongs to $\F^k$ for all
$k$, and also belongs to the limit flow $\F'$ if it exists.
Moreover, the mean curvature of a normalized blow-up sequence $\F^k$
is equal to 1 at $(0, 0)$.
\newline
(ii) A limit flow is an ancient flow.
\newline
(iii) The $\kappa$-noncollapsing is invariant under translation and
dilation of coordinates. Hence if $\F$ is $\kappa$-noncollapsing, so
is any blow-up sequence  with the same $\kappa$.
\newline
(iv) We always have $\kappa\le n$.
\newline
(v) It is easy to verify that $r_{x, t}$ depends continuously on
$(x, t)$ and $\F$. Namely if $\F^k$ converges locally smoothly to
$\F'$, and $\F^k\ni (x_k, t_k)\to (x_0, t_0)\in \F$, then we have
$r_{x_k, t_k}[\F^k]\to r_{x_0, t_0}[\F']$, where $r_{x_k,
t_k}[\F^k]$ is the radius in (3.2) relative to $\F^k$.
\newline
(vi) A semi-noncollapsing blow-up sequence may contain a collapsing
component, but we will rule out the case in \S5 (Remark 5.1). Also a
collapsing blow-up sequence may contain two or more collapsing
components.

\vskip5pt

\noo {\bf Notation}:
\newline
$\M$:\hskip10pt an oriented, compact $n$-dimensional differential
manifold without boundary.
\newline
$\F=\{F_t\}_{t\in [0, T)}$: \hskip10pt smooth solution to the mean
curvature flow with initial condition $F_0$, on a maximal time
interval $[0, T)$ for some $T\in (0, \infty]$.
\newline
$F_t=F(\cdot, t)$ ($0\le t<T$): \hskip10pt a smooth embedding of
$\M$ in $\R^{n+1}$.
\newline
$F_t$: \hskip10pt we also use $F_t$ to denote the hypersurface
$F_t(\M)$, and $F_0={F_t}_{|t=0}$.
\newline
$U_t$: \hskip10pt the open domain enclosed by $F_t(\M)$,
$U_0={U_t}_{|t=0}$, and  $\U=\{U_t\}_{t\in [0, T)}$.
\newline
$\bar U_t$: \hskip10pt closure of $U_t$, $\bar U_t=U_t\cup\p U_t$.
\newline
$\B_r(x)=\B^{(n+1)}_r(x)$: \hskip10pt the ball in $\R^{n+1}$ of
radius $r$ with center at $x$.
\newline
$\B^n_r(x')$: \hskip10pt the ball in $\R^n$ of radius $r$ with
center at $x'$.
\newline
$Q_r(x_0, t_0)$: \hskip10pt the parabolic cube $\B_r(x_0)\times
(t_0-r^2, t_0]$ in spacetime.
\newline
$|\cdot|_{\H^k}$:\hskip10pt  $k$-dimensional Hausdorff measure.

For a blow-up sequence $\F^k=\{F^k_t\}$ (or a limit flow
$\F'=\{F'_t\}$), we denote accordingly by $U^k_t$ (or $U'_t$) the
open set enclosed by $F^k_t$ (or $F'_t$).

In this paper, we use $x=(x_1, x_2 \cdots, x_{n+1})$ to denote a
point in $\R^{n+1}$, and use $x'=(x_1, x_2 \cdots, x_n)$ to denote a
point in $\R^{n}$.

\vskip0pt

\baselineskip=17.0pt

\section{ Convexity of limit flow}

\noo{\bf Lemma 4.1.}\ {\it If a blow-up sequence $\F^k$ converges
locally smoothly to a limit flow $\F'$, then $\F'$ is convex, namely
$U'_t$ is convex whenever it is nonempty.}

\noo{\bf Proof}. For any time $t_0$ such that $U'_{t_0}$ is not
empty, it suffices to prove that for any two interior points $z_0,
z_1\in U'_{t_0}$, the line segment $\overline{z_0z}_1$ lies in
$U'_{t_0}$.

By a proper translation and rotation of coordinates, we assume that
$z_0=(0, \cdots, 0, 1)$ and $z_1=(0, \cdots, 0, -1)$. Since $z_0,
z_1$ are interior points of $U_{t_0}'$, there exists small
$\delta>0$ such that $\B_{4\delta}(z_0), \B_{4\delta}(z_1)\subset
U'_{t_0}$. Denote
$$\Om=\{x\in\R^{n+1}:\
  \ {\Sigma}_{i=1}^nx_i^2<\delta^2(1+x^2_{n+1}), -1<x_{n+1}<1\}.$$
Since the initial hypersurface $F_0$ is smooth and mean convex, the
singularity set is strictly contained in the domain enclosed by
$F_0$. Hence $\Om\subset U^k_{t}$ for any large $k$, provided $-t$
is sufficiently large. If the line segment $\overline{z_0z}_1$ is
not contained in $U'_{t_0}$, we decrease $t$ to a moment $\bar
t=\bar t_k$ such that $\Om\subset U^k_{\bar t}$ but $\pom$ contacts
with the boundary of $U^k_{\bar t}$ at some point $z^*$, namely
 \begin{equation}
 \bar t=\sup\{s:\ \Om\subset U^k_{t}\ \forall\ t<s\}.\tag {4.1}
 \end{equation}
Since $\B_{4\delta}(z_0), \B_{4\delta}(z_1) \subset U'_{t_0}$, we
have $z^*\in \pom\cap\{-1<x_{n+1}<1\}$. Let $\lam$ and $\Lambda$ be
the least and largest principal curvatures of $F^k_{\bar t}$ at this
point. By comparing the principal curvatures of $\pom$ and
$F^k_{t_0}$,  we have $\lam\le -\frac 14\delta$ and $\Lambda\le
4\delta^{-1}$. Hence the mean curvature $H\le C\delta^{-1}$ by
(2.4). Applying these estimates to the curvature pinching (1.1) we
obtain $-\delta a_k\ge -\phi(C\delta^{-1}a_k)$. Note that $\delta>0$
is fixed, $a_k\to\infty$, and $\phi(t)/t\to 0$ as $t\to\infty$, we
reach a contradiction. $\hfill\square$

The above proof is inspired by an idea in [22] (pages 130-131,
[22]). We made some changes to make use of the curvature pinching
estimate (1.1). Note that the curvature pinching estimate in [11,
12] is written in a different way, but it is easy to see that it
implies (1.1).

\vskip0pt

\section{Smooth convergence}


We prove a local curvature estimate (5.3), from which it follows the
smooth convergence of a semi-noncollapsing blow-up sequence.

Given constants $\delta, d\in (0, 1]$, we say a mean convex flow
$\F$ has property $P^a_{\delta, d}$ at a point $p_0=(x_0, t_0)\in\F$
if (a) below is satisfied; and $\F$ has property $P^{a,b}_{\delta,
d}$ at $p_0$ if both (a) and (b) are satisfied; where
\newline
(a) $\exists\ y_0\in\R^{n+1}$ s.t. $|x_0-y_0|=d$, $\B_{\delta
d}(y_0)\subset U_{t_0}$, and the line segment $x_0y_0\subset \bar
U_{t_0}$;
\newline
(b) $F_t\cap \B_{10d} (x_0)\subset N_{\delta d/32}(\Gamma_t)$ for
some convex hypersurface $\Gamma_t$,  $\forall\ t\in (t_0-(\delta
d)^2, t_0]$, where $N_\delta$ denotes the $\delta$-neighborhood.

Assume $\F$ has property $P^{a,b}_{\delta, d}$ at $p_0$. Choose a
coordinate system such that $x_0$ is the origin and $y_0=(0, \cdots,
0, d)$. Then $\Gamma_t\cap\{|x'|<\frac{\delta d}{2}\}$ is the graph
of a Lipschitz continuous, convex function $\phi(\cdot, t)$ defined
in $\B^n_{\delta d/2}(0)$, $t\in (t_0-(\delta d)^2, t_0]$.

{\it Claim}:  $|\phi(x', t)|\le 8d$ for any $x'\in \B^n_{\delta
d/2}(0)$ and $t\in (t_0-(\delta d/8)^2, t_0]$.

Indeed, by assumption (a) we have $\phi(\cdot, t)\le d$ in
$\B^n_{\delta d/2}(0)$, for any $t\in (t_0-(\frac{\delta d}{8})^2,
t_0]$. Since the origin $0\in \F^k_{t_0}$, by the convexity of
$\phi$ we have $\phi(\cdot, t_0)\ge -2d$. For $t\in
(t_0-(\frac{\delta d}{8})^2, t_0]$, noting that $\F$ is a mean
convex flow, by assumption (b) we see that up to a small
perturbation, $\phi(\cdot, t)$ is increasing in $t$. If $\phi(\cdot,
t_0-(\frac{\delta d}{8})^2)\le -8d$ at some point in $\B^n_{\delta
d/2}(0)$, by convexity we have $\phi(\cdot,t_0-(\frac{\delta
d}{8})^2)\le -3d$ in $\B^n_{\delta d/4}(0)$. A simple application of
the comparison principle for the mean curvature flow implies that
$\phi(0,t_0)\le -d$, which is in contradiction with (b). The Claim
is proved.

Suppose $F_t\cap \B_{\delta d/2}(0)$ is also a graph of a function
$u(\cdot, t)$ for $t\in (t_0-(\frac{\delta d}{8})^2, t_0]$. By the
interior gradient estimate (2.5), we have
 \begin{equation}
 |D_{x'}u(x', t)|\le C'_{\delta, n} \tag {5.1}
 \end{equation}
for $|x'|<\delta d/4$ and $t\in (t_0-(\frac{\delta d}{16})^2, t_0]$.
(For application to a blow-up sequence below, we don't need to use
the interior gradient estimate (2.5) to get (5.1). Indeed, by Lemma
5.1 and the smooth convergence in the proof of Lemma 5.2, the
estimate (5.1) follows from the Lipschitz continuity of the convex
function $\phi$). Hence by the interior estimate (2.6), we have
 \begin{equation}
 \Lambda_{\F}(0) \le C''_{\delta,n}/d,\tag {5.2}
 \end{equation}
where $\Lambda_{\F}(p)=|\lam_1(p)|+\cdots+|\lam_n(p)|$, and
$\lam_1(p), \cdots, \lam_n(p)$ are the principal curvatures of $F_t$
at $x$, $p=(x, t)\in\F$. The constants $C'_{\delta, n}$ and
$C''_{\delta, n}$ depend only on $\delta$ and $n$ but are
independent of $d$, which can be easily seen by making the dilation
$x\to x/d$.

\vskip5pt

Next we prove a technical lemma. Denote $\C_{x_0, \delta}=\{x\in
\R^{n+1}:\ x=x_0+t(y-x_0), \ t\in (0, 1), \ y\in \B_\delta(0)\}$ an
open round cone with vertex at $x_0$.

\vskip5pt

\noo{\bf Lemma 5.1.}\ {\it Let $U$ be an open set in $\R^{n+1}$ with
smooth boundary $F$. Suppose $\B_r(0)\subset U\subset \B_R(0)$. Then
$\exists\ \eps>0$, depending only on $n,r, R$, such that the
following conclusion holds: If at any point $x_0\in F$ with
$\C_{x_0, r/4}\subset U$, the principal curvatures of $F$ $\ge
-\eps$, then $\C_{x_0, r/2}\subset U$ for all $x_0\in F$. }

\noo{\bf Proof}. If the lemma is not true, assume $\inf\{|x|:\ x\in
F:\ \C_{x, r/2}\not\subset U\}$ is attained at $x_0$. Let $P$ be the
two-dimensional plane determined by $0$, $x_0$, and $\nu$, where
$\nu$ is the normal of $F$ at $x_0$. Let $\ell$ be the intersection
of $P$ with $F$. Then the tangent line of $\ell$ at $x_0$ is also
tangential to $\p \B_{r/2}(0)$, and a tangent plane of $F$ at $x_0$
is also a tangent plane of $\p \B_{r/2}(0)$. Let us parameterize
$\ell=\{x(s):\ s\ge 0\}$ by its arclength, starting at $x_0$, and
let $\alpha(s)=|x(s)|$. Then when $s>0$ small, $\alpha$ is
decreasing in $s$. Let $[0, s_0]$ be the maximal interval in which
$\alpha(s)\le \alpha(0)$.  By our choice of $x_0$, we have
$\C_{x(s), r/2}\subset U$ for all $s\in (0, s_0]$. Hence by
assumption, the principal curvatures of $F$ at $x(s)$ $\ge -\eps$.
Since $\C_{x(s), r/2}\subset U$, we have $\lan \nu(s), -x(s)\ran\ge
c_0>0$, where $\nu(s)$ is the normal of $F$ at $x(s)$. Hence the
curvature of the curve $\ell>-C\eps$. When $\eps>0$ is sufficiently
small, the curve $\ell$ must intersect with the ball $\B_r(0)$ (as
the tangent line of $\ell$ at $x_0$ is tangential to $\p
\B_{r/2}(0)$), which is in contradiction with the assumption
$U\supset \B_r(0)$. \hfill $\square$

\vskip5pt

When the condition $U\subset \B_R(0)$ is dropped, we consider the
set $U\cap\B_R(0)$. By approximation we have

\vskip5pt

\noo{\bf Lemma 5.1$'$.}\ {\it Let $U$ be an open set in $\R^{n+1}$
with smooth boundary $F$. Suppose $\B_r(0)\subset U$. Then $\exists\
\eps>0$, depending only on $n,r, R$, such that the following
conclusion holds: If at any point $x_0\in F$ with $|x_0|<R$ and
$\C_{x_0, r/4}\subset U$, the principal curvatures of $F$ $\ge
-\eps$, then $\C_{x_0, r/2}\subset U$ for any $x_0\in F$ with
$|x_0|<R$. }

\vskip5pt

Lemma 5.1 (and 5.1$'$) implies that for any given $R>0$, there is a
component of $F\cap \B_R(0)$ which can be represented as a Lipschitz
continuous radial graph over the unit sphere, and there is no other
components inside this component. We call this component the {\it
noncollapsing component} of $\F$.

\vskip5pt

Let $\F^k$ be a blow-up sequence. Assume $U^k_0$ contains a common
ball $\B_r(y_0)$ for all large $k$. We fix a small $\delta>0$ such
that $\F^k$ has the property $P^a_{\delta, d}$ at the origin for
some $d\le 1$. Note that for any smooth mean curvature flow $\F$ and
any $p\in\F$, $\F$ has property $P^a_{\delta, d}$ at $p$ with
$\delta\to 1$ as $d\to 0$. Hence for any $p\in\F^k$, we can assign a
constant $d=d_{k,p}>0$ such that $\F^k$ has the property
$P^a_{\delta, d}$ at $p$. We have freedom in choosing the value of
$d_{k, p}$ but we can choose the largest $d\in (0, 1]$ such that
$\F^k$ has the property $P^a_{\delta, d}$ at $p$.  Then for fixed
$k$, $d_{k, p}$ is upper semicontinuous in $p$, namely $d_{k,
p_0}\ge \lim_{p\to p_0} d_{k,p}$.

\vskip5pt

\noo{\bf Lemma 5.2.}\ {\it At any point $p=(x, t)\in\F^k$, the
principal curvatures of $F^k_t$ at $x$ satisfy
 \begin{equation}
 \Lambda_{\F^k}(x,t) < 2 C''_{\delta,n}/d_{k,p},\tag {5.3}
 \end{equation}
provided $k$ is sufficiently large. }

\vskip5pt

\noo{\bf Proof}. If the lemma is not true,  let
$$\tilde t_k=\sup\{t':\ \ \Lambda_{\F^k}(p)\le 2 C''_{\delta,n}/d_{k,p},
                 \ p=(x,t)\in\F^k,\ \forall\ x\in F^k_t, t\le t'\}.$$
For any given $k$, since $d_{k,p}$ is upper semi-continuous in $p$
and $\F^k$ is a smooth mean curvature flow, there exists a point
$\tilde p_k=(\tilde x_k, \tilde t_k)\in \F^k$ such that
 \begin{align}
 & \Lambda_{\F^k}(\tilde p_k)\ge 2C''_{\delta,n}/d_{k, \tilde p_k},
  \tag {5.4}\\
 & \Lambda_{\F^k}(p)\le 2C''_{\delta,n}/d_{k,p}
 \ \ \forall\ p=(x, t)\in\F^k\ \text{with}\ t < \tilde t_k.\tag
 {5.5}
 \end{align}
Let $\tilde \F^k=\D_{\tilde a_k, \tilde p_k}(\F^k)$, where $\tilde
a_k=d_{k, \tilde p_k}^{-1}$ and $\D$ is the dilation in (3.1). Then
for any point $p=(x, t)\in \tilde \F^k$ with $t< 0$, the principal
curvatures of $\tilde \F^k$ satisfy $ \Lambda_{\tilde\F^k}(x,t)\le 2
C''_{\delta,n}/\tilde d_{k,p}$,  where $\tilde d_{k,p}\in (0, 1]$ is
the largest constant such that $\tilde \F^k$ has the property
$P^a_{\delta, \tilde d_{k,p}}$ at $p$. By our definition, $\tilde
d_{k,p}=1$ at $p=(0,0)$. Hence there is a ball $\B_\delta(\tilde y)$
contained in the interior of $\tilde F^k_t$ for all $t\le 0$.

Now we apply Lemma 5.1 to $\tilde F^k_t$ with the ball $\B_r(0)$ in
Lemma 5.1 replaced by $\B_{\delta/8}(\tilde y)$. Note that the
condition $P^a_{\delta, d}$ is weaker than the cone condition in
Lemma 5.1. By (5.5) and the curvature pinching estimate (1.1), at
any point $p=(x,t)\in \tilde \F^k$ with $t<0$,  the principal
curvatures of $\tilde\F^k$ $\ge -\eps$ with $\eps\to 0$ as
$k\to\infty$. Hence by Lemma 5.1, $\tilde F^k_t\cap\B_R(0)$ is a
Lipschitz continuous radial graph (with respect to $\tilde y$) when
$t<0$. Therefore $\tilde \F^k$ sub-converges locally smoothly to a
limit flow $\tilde \F'$, which is convex by Lemma 4.1. Hence when
$k$ is sufficiently large, $\tilde F^k_t$ satisfies the conditions
(a) and (b) above. Sending $t\nearrow 0$ we obtain by (5.2) the
curvature estimate $\Lambda_{\tilde \F^k}(0)\le C''_{\delta,n}$,
which is in contradiction with (5.4). \hfill$\square$

\vskip5pt

From estimate (5.3) and Lemma 5.1, we see that $F^k_t\cap \B_R(0)$
is a Lipschitz continuous radial graph (with respect to some point)
for $t\in (-R^2, 0]$. Therefore we obtain the following local smooth
convergence.

\vskip5pt

\noo{\bf Corollary 5.1.}\ {\it If a blow-up sequence $\F^k$ is
semi-collapsing, then the noncollapsing component of $\F^k$
sub-converges locally smoothly to a convex limit flow. }

\vskip5pt

\noo{\bf Remark 5.1}. To conclude that $\F^k$ itself converges
locally smoothly, we need to rule out the possibility that $\F^k$
contains two or more components, one is the semi-noncollapsing
component and the others are collapsing (the argument in \S4 implies
that a blow-up sequence cannot contain two or more
semi-noncollapsing components).

Let $p_0=(x_0, t_0)$ be a point inside the collapsing component (we
need to translate $\F^k$ slightly such that the collapsing component
of $\F^k$ share a common point $p_0$ for all $k$), and let
$\B_{2r}(p_1)$ be a ball in $U'_0$ (we assume without loss of
generality that $U'_0\ne\emptyset$). Let $\C=\{q=tp_0+(1-t)p:\ t\in
(0, 1),\ p\in \B_r(p_1)\}$ be a convex cone with vertex at $p_0$. As
in the proof of Lemma 4.1, let
$$\hat t_k=\sup\{s:\ \C\subset U^k_t\ \forall\ t<s\}.$$
At time $\hat t_k$, there is a point $\hat x_k\in\p\C\cap F^k_{\hat
t_k}$, $\hat x_k\ne x_0$. Let $P_k$ be the tangent plane of
$F^k_{\hat t_k}$ at $\hat x_k$. Assume $P_k\to P_0$ as $k\to\infty$.
Making the translation $(x, t)\to (x, t-\hat t_k)$, we get a new
sequence $\hat\F^k$. Obviously $\hat\F^k$ is semi-noncollapsing. By
Corollary 5.1, it converges locally smoothly to a convex limit flow
$\hat \F'$. The smooth convergence implies that $P_0$ is the tangent
plane of $\hat F'_0$ at $0$, where $\hat F'_0$ is the hypersurface
of $\hat\F'$ at time $t=0$. The convexity implies the interior $\hat
U'_0$ of $\hat F'_0$ lies on one side of $\hat P_0$. But on the
other hand, the ball $\B_{r}(p_1)$ is tangent to $P_0$ and
$\B_{2r}(p_1)$ is contained in $\hat U'_0$, we reach a
contradiction.

Therefore we have proved the following theorem. We point out that
Corollary 5.1 is sufficient for the treatment in the next section to
rule out collapsing.

\vskip5pt

\noo{\bf Theorem 5.1.}\ {\it Let $\F^k$ be a semi-noncollapsing
blow-up sequence. Then it converges locally smoothly to a convex
limit flow. In particular, $\{F'_t\}=\{\p U'_t\}$ is a convex
solution to the mean curvature flow. }

\vskip5pt

\noo{\bf Remark 5.2.}\ The arguments above also applies to a
collapsing blow-up sequence $\F^k$ provided $\F^k$ is locally
separate. More precisely, let $\F^k=\{F^k_t\}_{t\in [-\hat T, 0)}$
be a blow-up sequence (here we restrict the time $t$ to a fixed
interval). Denote
 \begin{equation}
 \C_{R}=\{x=(x', x_{n+1})\in\R^{n+1}:\ |x'|<R\} \tag {5.6}
 \end{equation}
a cylinder in $\R^{n+1}$. Suppose for all $t\in [-\hat T, 0]$,
$$F^k_t\cap \C_{R} \subset \{|x_{n+1}|\le 10^{-3}\},$$
and $F^k_t\cap \C_{R}$ consists of two disconnected components,
$$F^k_t\cap \C_{R}=G^k_{1, t}\cup G^k_{2,t},$$
such that the projection of both $G^k_{1, t}$ and $G^k_{1, t}$ in
the plane $\{x_{n+1}=0\}$ covers the ball $\B^n_R(0)$. Then $G^k_{1,
t}$ divides the ball $\B_{R+t}(0)$ into two separate open sets (we
choose the variable radius $R+t$ such that the sets $U^k_{1,t}$ and
$\B^+_k$ depend continuous on $t$. The continuous dependence is
needed in the proof of Lemma 5.2). Let $\B^+_k$ be one of the open
sets which contains $U^k_t\cap\B_{R+t}(0)$, and let
$U^k_{1,t}=U^k_t\cup \B^+_k$. The proof of Lemma 4.1 implies that
the limit of $\B^+_k$ is convex. To extend Lemma 5.2 to the current
case, we define $d_{k, p}$ using the open set $U^k_{1,t}$ (for $p\in
G^k_{1, t}\cap \B_{R+t}(0)$, for other points on $\F^k$, we define
$d_{k,p}$ as in Lemma 5.2). Then the argument above implies that
$G^k_{1, t}$ converges locally smoothly to a convex hypersurface.
Similarly $G^k_{2, t}$ converges locally smoothly to a convex
hypersurface.

\vskip0pt

\section{Semi-noncollapsing}

In this section we prove that every blow-up sequence is
semi-noncollapsing. The proof consists of two steps, in the first
one we show that if every collapsing blow-up sequence is double
sheeting, then there is no collapsing blow-up sequence (Lemma 6.2).
In the second one we show that every collapsing blow-up sequence
must be double sheeting (Lemma 6.4).

\vskip5pt

First we introduce the notion of {\it double sheeting}. We say a
blow-up sequence $\F^k$ is double sheeting if for any $R>0$,
$\eps>0$, there exists $k_{R,\eps}>1$ such that for all $k\ge
k_{R,\eps}$, $t\in (-\eps^{-1}, -\eps)$,
 \begin{equation}
F^k_t\cap \C_{R}=G_{f^k}\cup G_{g^k},\tag {6.1}
 \end{equation}
where $G_{f^k}$ and $G_{g^k}$ are graphs of functions $f^k=f^k(x',
t)$ and $g^k=g^k(x', t)$, in the form
\begin{align}
G_{f^k} & =\{x=(x', x_{n+1}):\ x_{n+1}=f^k(x', t),\ |x'|<R\}, \tag {6.2}\\
G_{g^k} & =\{x=(x', x_{n+1}):\ x_{n+1}=g^k(x', t),\ |x'|<R\},
\nonumber
\end{align}
such that
 \begin{equation*}
 U^k_t\cap\C_{R}=\{g^k(x', t)<x_{n+1}<f^k(x', t)\}.
 \end{equation*}
Moreover, for any $|x'|<R$, $-\eps^{-1}<t<-\eps$, and $k\ge k_{R,
\eps}$,
\begin{equation}
{\sum}_{i=0}^4\big\{|D_{x'}^if^k(x', t)|+|D^i_{x'}g^k(x', t)|\big\}
 \le \eps. \tag{6.3}
 \end{equation}
By the mean convexity of $F^k_t$,  $f^k$ is decreasing, and $g^k$ is
increasing, in $t$.

\vskip5pt

Let $\F^k=\{F^k_t\}$ be a blow-up sequence. For any $x\in \bar
U^k_t$ (where $\overline U$ is the closure of $U$), denote
\begin{equation}
 r_{k,x,t}=
 \sup\{\rho:\ x\in \overline \B_{\rho}(z)
 \ \text{and}\ \B_\rho(z)\subset U^k_t
 \ \text{for some}\ z\in U^k_t\}.\tag {6.4}
   \end{equation}
Recall that by our definition in \S3, the origin in space-time is a
point in $\F^k$. For any given $\tau<0$, if there is a positive
constant $c_0$ such that
$$r_{k, 0, \tau}\ge c_0\ \ \forall\ k$$
by Theorem 5.1, $\F^k$ converges locally smoothly to a convex limit
flow. If
$$r_{k, 0, \tau}\to 0\ \ \text{as}\ k\to\infty, $$
let $\wtt \F^k=\D_{a_k, p_k}(\F^k)$, where $a_k=1/r_{k, 0, \tau}$,
$p_k=(0, \tau)$, and $\D$ is the dilation in (3.1). By Theorem 5.1,
$\wtt \F^k$ converges locally smoothly to a convex limit flow
$\wtt\F'=\{\wtt F'_t\}_{t\in\R^1}$.

\vskip5pt

\noo{\bf Lemma 6.1}. {\it The limit flow $\wtt\F'$ is a pair of
parallel planes, namely, $\wtt F'_t$ is a pair of parallel planes
for all $t\in\R^1$.}

\noo{\bf Proof}. By (6.4) and the smooth convergence of $\wtt\F^k$
to $\wtt\F'$ , $\wtt U'_0$ does not contain a ball $\B_r(z)$ such
that $r>1$ and $0\in \overline \B_r(z)$. Hence $\wtt U'_0$ cannot be
a half space and $\wtt F'_0$ cannot be a (single) plane.

If the limit flow $\wtt\F'$ is not a pair of parallel planes, then
the mean curvature of $\wtt \F'$ is positive everywhere and it will
pass through the origin in finite time. But on the other hand, since
$\tau<0$, by the parabolic dilation above, the origin is contained
in $U'_t$ for all $t>0$. This is a contradiction. \hfill$\square$

\vskip5pt

Note that by Hamilton's maximum principle for tensors, the second
fundamental form of a convex mean curvature flow has constant rank
[12]. It implies the mean curvature is positive everywhere.

\vskip5pt

Therefore if $r_{k, 0, \tau}\to 0$ for some $\tau<0$, then $\F^k$ is
a collapsing blow-up sequence, and for any given $t<0$, $U^k_t$ does
not contain a common ball for all large $k$ (this property will be
used in the argument below). Indeed, if $U^k_t$ (for some $t<0$)
contains a common ball for a subsequence $k\to\infty$, by Theorem
5.1, $\F^k$ converges locally smoothly to a convex limit flow $\F'$.
Suppose $\F'$ becomes extinct at $T'$. By the smooth convergence
$\F^k\to\F'$ and the convexity of $\F'$, it is not hard to show that
$\F^k\cap\B_R(0)$ becomes extinct before $T'+\eps_k$ with $\eps_k\to
0$ as $k\to\infty$.

\vskip5pt

\noo{\bf Lemma 6.2}. {\it  If every collapsing blow-up sequence is
double sheeting, then collapsing blow-up sequence does not occur. }

\vskip5pt

\noo{\bf Proof}. Let $\F=\{F_t\}_{t\in [0, T)}$ be a mean convex
flow which develops first time singularity at time $T$. Let
$\U_T=\U\cap \{t\le T\}$ (see notation in \S3). Regard $\U_T$ as a
domain in the space-time. For any small $r>0$, and any given point
$p=(x_0, t_0)\in \p\U_T$ with $t_0\ge\frac 12 T$, denote
 \begin{align*}
 w_p(r) &=
 \sup\{\rho:\ x_0\in \overline\B_\rho(z)
    \ \text{and}\ \B_\rho(z)\subset U_{t_0-r^2}\},\\
 w(r) &=\inf_p w_p (r).
 \end{align*}
Since $F_t$ is smooth, closed, and mean convex for $t<T$, we have
\newline
(i) $w$ is continuous and non-decreasing in $r$;
\newline
(ii) $w(r)>0$ for any $r>0$;
\newline
(iii) $w(r)$ is attained at some point $(x_r, t_r)$.

If there is a collapsing blow-up sequence, we have
 \begin{equation}
 {\underline{\lim}}_{r\to 0}\frac {w(r)}{r}=0.\tag {6.5}
 \end{equation}
Thus there must be a sequence $r_i\to 0$ such that $\frac
{w(r_i)}{r_i}\to 0$ and
\begin{equation*}
 \frac {w(r_i)}{r_i}\le 2\frac{w(3r_i)}{3r_i}.
\end{equation*}

Suppose $w(r_i)$ is attained at $p_i=(x_i, t_i)$. By definition we
have $w_{p_i}(r_i)=w(r_i)$ and $w_{p_i}(3r_i)\ge w(3r_i)$. Hence
\begin{equation}
 \frac {w_{p_i}(r_i)}{r_i}\le 2\frac{w_{p_i}(3r_i)}{3r_i}. \tag {6.6}
\end{equation}

Let $\F^k=\D_{r_i^{-1}, p_{i}}(\F)$ be the blow-up sequence of $\F$
at $p_{i}$. Since $\frac{w_{p_i}(r_i)}{r_i}\to 0$, $r_{k, 0,
\tau}\to 0$ at $\tau=-1$, where $r_{k, 0, \tau}$ is the number given
in (6.4), relative to the blow-up sequence $\F^k$. Hence $\F^k$ is
collapsing. Hence by assumption, $\F^k$ is double sheeting. Suppose
the two sheets are given by (6.2). Denote
 \begin{align*}
  u^k(x, t)&=(f^k-g^k)(x, t),\\
 \hat u^k(x, t)& =u^k(x, t)/u^k(0, -1),
 \end{align*}
Then $u^k$, together with its first and second derivatives,
converges locally uniformly to 0 in the parabolic cube $Q_R(0, -1)$,
for any $R>1$. Note that $f^k$ and $g^k$ satisfy the parabolic
equation (2.2). Hence $u^k$ and $\hat u^k$ satisfy a parabolic
equation in $Q_R(0, -1)$ which converges as $k\to \infty$ locally
uniformly to the heat equation. Since $\hat u^k(\cdot,t)$ is
decreasing in $t$, the Harnack inequality
  \begin{equation}
  \sup_{Q_R(0, -1)} \hat u^k\le C\inf_{Q_R(0, -1)} \hat u^k
  \tag {6.7}
 \end{equation}
holds, where the constant $C$ is independent of both $k$ and $R$,
provided $k$ is sufficiently large such that $\hat u^k$ is
well-defined in $Q_{2R}(0, -1)$ and satisfies the double sheeting
condition. Hence by the regularity of parabolic equation, $\hat u^k$
converges locally smoothly to a function $\hat u$ which is
nonincreasing in $t$, and satisfies the heat equation in $\R^n\times
(-\infty, -1]$. Applying the above Harnack inequality to $\hat
u-\inf_{\R^n} \hat u$ in $Q_R(0, -1)$ and sending $R\to\infty$, we
concludes that $\hat u -\inf \hat u\equiv 0$, namely $\hat u\equiv
\inf \hat u$ in $\R^n$. Hence $\hat u^k$ converges locally smoothly
to $\hat u(0,-1)=1$ in $\R^n$. On the other hand, by (6.6) we have
$w_{p_i}(3r_i)\ge \frac 32 w_{p_i}(r_i)$. Hence
$$\hat u(0, -9)=\lim_{k\to\infty} \hat u^k(0, -9)
 \ge \frac 32\lim_{k\to\infty}\hat u^k(0, -1)=\frac 32,$$
we reach a contradiction. $\hfill\square$

The above proof is inspired by the proof of Theorem 9.2 in [21]
(pages 688-690). Next we show that every collapsing blow-up sequence
must be double sheeting. First we prove

\noo{\bf Lemma 6.3.} {\it  Let $\F^k$ be a collapsing blow-up
sequence. Then for any $\eps>0$, there exists $k_\eps>1$ such that
after a rotation of axes,
 \begin{equation}
 F^k_t\cap\C_{\eps^{-1}}\subset\{|x_{n+1}|<\eps\}
 \tag{6.8}
 \end{equation}}
for all  $t>-\eps^{-1}$ and $k>k_\eps$.

 \noo {\bf Proof}. For any $m>1$, let
  \begin{equation*}
  t_{k,m}=\sup\{\tau:\ r_{k, 0, \tau}>1/m\},
   \end{equation*}
where $r_{k, 0, \tau}$ is given in (6.4). Then $t_{k, m}\to-\infty$.
Translating $\F^k$ forward in time by $t_{k,m}$, such that $(0,
t_{k,m})$ becomes the origin in space-time,  we get a blow-up
sequence $\F^k_m$ such that $\{0\}\in\B_{1/m}(z)\subset (U^k_m)_0$
(the interior of $\F^k_m$ at $t=0$). By Theorem 5.1, $\F^k_m$
converges locally smoothly to a convex limit flow $\F'_m$.

From the proof of Lemma 6.1, $\F'_m$ is a pair of parallel planes,
namely $(F'_m)_t$ is a pair of parallel planes for all $t\in\R^1$.
For if not, the flow $\F'_m$ would pass through the origin in finite
time. But on the other hand, we have $t_{k, m}\to-\infty$, which
means the origin is contained in $(U'_m)_t$ for all $t>0$.

By a rotation of the axes, we assume the parallel planes is given by
$\{x_{n+1}=c_1\}$ and $\{x_{n+1}=c_2\}$, for some constants $c_1,
c_2$ satisfying $c_1\le 0\le c_2$ and $c_2-c_1\le 5/m$. Hence for
any given $R$,
$$\C_{R}\cap F^k_{t_{k, m}}\subset \{|x_{n+1}|<10/m\} $$
provided $k$ is sufficiently large. By the mean convexity, we have
$$ \C_{R}\cap F^k_t\subset \{|x_{n+1}|<10/m\}
 \ \ \ \forall\ t>t_{k,m} $$
provided $k$ is sufficiently large. We obtain (6.8) by choosing
$m=10/\eps$ and $R=\eps^{-1}$. \hfill$\square$

\vskip5pt

For any $m>1$, $\F^k_m$ is semi-noncollapsing. The above proof
implies that $F^k_t$ is asymptotically double sheeting when $t\in
(t_{k, m}, t_{k, 2m})$, for any sufficiently large $m$. Therefore in
the following it suffices to consider the time $t$ in the range
$0>t>t_{k, m}$.

For any $\tau <0$, let $\wtt \F^k=\D_{a_k, p_k}(\F^k)$, where
$a_k=1/r_{k, 0, \tau}$, $p_k=(0, \tau)$. By Lemma 6.1,  $\wtt \F^k$
converges locally smoothly to a pair of parallel planes, given by
$\{x_{n+1}=c_1\}$ and $\{x_{n+1}=c_2\}$. Hence for any $R>0$, when
$k$ is sufficiently large, $\wtt F^k_t\cap\C_{R}$ consists of two
separate parts, which are graphs of $\wtt f^k$ and $\wtt g^k$, given
by $x_{n+1}=\wtt f^k(x', t)$ and $x_{n+1}=\wtt g^k(x', t)$, and
$\wtt f^k\to c_1$, $\wtt g^k\to c_2$ uniformly for $x'\in \B_R(0)$
and $-R^2<t<0$.

Rescaling back we see that for any time $t<0$, $F^k_t\cap\C_{r}$ is
double sheeting locally for small $r>0$, namely it consists of two
separate parts $G_{f^k}$ and $G_{g^k}$, as given in (6.2). Denote
$$V_{k,\tau, r}=\{x\in U^k_{\tau}:\ \ d_{0, x}[U^k_{\tau}]<r\},$$
\vskip-15pt

\noo where
 \begin{equation}
 d_{x, y}[U]=\inf\{ \text{arclength of paths contained in}\ U,
  \text{connecting}\ x\ \text{to}\ y\}. \tag {6.10}
  \end{equation}
From the above discussion, for $r>0$ small, $\big(\p V_{k,\tau,
r}\big)\cap F^k_t$ consists of two separate components
$\mu^+_{k,\tau,r}$ and $\mu^-_{k,\tau,r}$.  Let
 \begin{align*}
 \rho_{k,\tau}=\sup\{\bar r:\ \ &
 \text{the two separate parts $\mu^+_{k,\tau,r}$ and
 $\mu^-_{k,\tau,r}$}\\
  & \text{keep separate for all}\ r\in (0, \bar r)\}.
  \end{align*}
We point out that $\mu^+_{k,\tau,r}$ and $\mu^-_{k,\tau,r}$ are the
graphs $G_{f^k}$ and $G_{g^k}$ when $r$ is small, but may not be
graph when $r$ is large.  Let $y_{k,\tau}$ be the point at which the
two parts $\mu^+_{k,\tau,\rho_{k,\tau}}$ and
$\mu^-_{k,\tau,\rho_{k,\tau}}$ meet.

\vskip5pt

\noo{\bf Lemma 6.4}. {\it If $\F^k$ is a collapsing blow-up
sequence, then it is double sheeting. }

\noo{\bf Proof}. We need to consider three different cases.

\noo {\it Case 1}: $|y_{k,\tau}|\to\infty$ for all $\tau<0$. In this
case, for any given $R>0$, $\C_{R}\cap F^k_t$ consists of two
separate components for all $t\le \tau$. By Remark 5.2 and Lemma
6.3, each of the two components converges locally smoothly to the
plane $\{x_{n+1}=0\}$. Hence $\F^k$ is double sheeting.

\vskip5pt

\noo {\it Case 2}: $|y_{k,\tau}|\to a_0>0$ for some $\tau<0$.  As in
Case 1, $\C_{R}\cap F^k_t$ (for any $R<a_0$) consists of two
separate components which converges locally smoothly to the plane
$\{x_{n+1}=0\}$. Let $\ell_{k,\tau}$ be the shortest path contained
in $\bar U^k_\tau$ which connects $0$ to $y_{k,\tau}$. Then by
definition, the arclength of $\ell_{k, \tau}$ is equal to $\rho_{k,
\tau}$. Since $F^k_t$ is smooth for all $t<0$, there is a point
$\hat y_{k,\tau}\in\ell_{k, \tau}$ such that
 \begin{equation}
 r_{k, \hat y_{k,\tau}, \tau}
  = \frac {1}{10} d_{\hat y_{k,\tau},y_{k,\tau}}[U^k_\tau].\tag{6.11}
  \end{equation}
Let $\hat \F^k=\D_{a_k, p_k}(\F^k)$, where $a_k=1/r_{k, \hat
y_{k,\tau}, \tau}$, $p_k=(\hat y_{k,\tau}, \tau)$.  We get a new
blow-up sequence which converges to a convex limit flow $\hat\F'$
with a nonempty interior. By (6.11), the mean curvature of $\hat
F'_0$ is strictly positive at the point $D_{a_k, p_k}(y_{k,\tau})$,
i.e. $H\ge C_0>0$. Scaling back, we see that the mean curvature of
$F^k_\tau$ at $y_{k, \tau}$ is greater than $C_0a_k\to\infty$. Hence
$$\frac{\p}{\p t} \rho_{k, t}\le -Ca_k\to-\infty$$
as $k\to\infty$. By Lemma 6.5 below, $\rho_{k, t}\le 4|y_{k,t}|$.
Hence $\rho_{k, t}\to 0$ in very short time, which means the blow-up
sequence $\F^k$ passes through the origin in very short time, and so
the origin $0\not\in F^k_0$, we reach a contradiction.

\vskip5pt

\noo {\it Case 3}: If there exists a $\tau<0$ such that
$|y_{k,\tau}|\to 0$, we make a parabolic dilation such that
$|y_{k,\tau}|=1$. The new blow-up sequence is still collapsing and
so we are in the same situation as in Case 2, which implies the
blow-up sequence $\F^k$ passes through the origin immediately.
$\hfill\square$

\vskip5pt

By our choice of the point $y_{k, \tau}$, there are at least $n-1$
vanishing principal curvatures of $\hat F'_0$ (the hypersurface at
$t=0$ of the limit flow $\hat\F'$ in Case 2) at the limit point of
$D_{a_k, p_k} (y_{k,\tau})$. Therefore by the constant rank of the
second fundamental form of $\hat \F'$ [12], the limit flow $\hat
\F'$ is the product of $\R^{n-1}$ with a convex solution to the
curvature shortening flow. By (6.11), $C_0$ is an absolute constant.

\vskip10pt

\noo{\bf Lemma 6.5}. {\it Let $\rho_{k, \tau}$ and $y_{k,\tau}$ be
as above. If $|y_{k, \tau}|$ is uniformly bounded, then we have
$\rho_{k, \tau}\le 4|y_{k,\tau}|$.}

\noo{\bf  Proof.}  Let $\ell_{k,\tau}$ be the shortest path in
$U^k_\tau$ connecting $0$ to $y_{k,\tau}$. For a sequence of points
$z_k\in\ell_{k,\tau}$, let $r_k=r_{k, z_k, \tau}$ be given in (6.4)
and let $d_k=d_{z_k, y_{k, \tau}}$ be the arelength of $\ell_{k,
\tau}$ from $z_k$ to $y_{k,\tau}$.

Let $z_0=0$, we choose a sequence of points $z_m=z_{m,k}\in
\ell_{k,\tau}$ for $m=1,2,\cdots$ such that the angle between
$\xi(z_m)$ and $\xi(z_{m-1})$ is equal to $\frac\pi {6}$, and for
any point $z\in \ell_{k,\tau}$ between $z_{m-1}$ and $z_m$, the
angle between $\xi(z)$ and $\xi(z_{m-1})$ is less than or equal to
$\frac\pi {6}$, where $\xi(z)$ denote the tangent vector of
$\ell_{k, \tau}$. If there is no point $z\in \ell_{k,\tau}$ such
that the angle between $\xi(z)$ and $\xi(z_0)$ $\ge \frac\pi {6}$,
Lemma 6.5 is obvious.

{\it Claim}: $d_{z_{m}, z_{m+1}}\le \frac 14 d_{z_{m-1}, z_m}$
uniformly in $m$, provided $k$ is sufficiently large.

Indeed, if there is a sequence of $m_k$ such that $d_{z_{m_k, k},
z_{m_k+1, k}}\ge \frac 14 d_{z_{m_k,k}, z_{m_k-1, k}}$ (recall that
$z_m=z_{m, k}$), by a parabolic dilation, we may assume that
$d_{z_{m_k,k}, z_{m_k-1, k}}=1$.  If $\F^k$ is collapsing after the
parabolic dilation, then $\F^k\cap \B_{1/8}(z_{m_k,k})$ consists of
two separate components and both converges smoothly to a hyperplane
(here $\B_{1/8}$ is the ball after the dilation). If $\F^k$ is
semi-noncollapsing after the parabolic dilation, by Theorem 5.1 it
converges locally smoothly to a convex solution. In both case we
reach a contradiction, as $\ell_{k,\tau}$ is the shortest path
connecting $z_{m-1}$ to $z_{m+1}$ but $\text{osc}_{z\in
\ell_{z_{m-1}, z_{m+1}}} \xi(z)\ge \frac \pi 6$, where
$\ell_{z_{m-1}, z_{m+1}}$ is the part of $\ell_{k,\tau}$ between
$z_{m-1}$ and $z_{m+1}$.

Hence we obtain
$$\rho_{k, \tau}={\sum}_{m\ge 1} d_{z_{m-1}, z_m}[U^k_\tau]
   <2d_{z_0, z_1}[U^k_\tau]. $$
Note that for any $R<|y_{k,\tau}|$, $\F^k\cap \B_{1/8}(0)$ consists
of two separate components, and each converges  smoothly to a
hyperplane. Hence we have $|y_{k, \tau}|\ge \frac 12d_{z_0, z_1}$.
\hfill$\square$

\vskip5pt

We have thus proved that every blow-up sequence is
semi-noncollapsing.

 \vskip10pt

\noo{\bf Corollary 6.1}. {\it  A limit flow cannot be a pair of
parallel planes. Moreover, the grim reaper is not a limit flow.}

\noo{\bf Proof}. If a blow-up sequence $\F^k$ converges to a
multiplicity two plane, then $\F^k$ is collapsing and must be double
sheeting, which was ruled out in Lemma 6.2. If $\F^k$ converges to a
pair of parallel planes, we can choose a sequence $a_k\to 0$ such
that $\D_{a_k, 0}(\F^k)$ is a collapsing blow-up sequence. But any
blow-up sequence is semi-noncollapsing, we also get a contradiction.

If the grim reaper is a limit flow, then by a proper translation,
there is a limit flow which is a pair of parallel planes.
$\hfill\square$

\vskip10pt

\baselineskip=16.5pt

\section{Proof of Theorems 1.1-1.3}

Theorem 1.1 follows from Theorem 5.1 and the semi-noncollapsing in
\S6.

\vskip5pt

\noo{\bf Proof of Theorem 1.2}.  If there is a sequence of points
$p_k=(x_k, t_k)$ such that
 \begin{equation}
 r_{x_k, t_k}H(x_k, t_k)= \delta_k\to 0 ,\tag {7.1}
 \end{equation}
where $r_{x,t}$ is defined in (3.2), we consider the blow-up
sequence $\F^k=\D_{a_k, p_k}(\F)$, where $a_k=1/r_{x_k, t_k}$. By
Theorem 5.1, $\F^k$ converges locally smoothly to a convex limit
flow $\F'$. By (7.1), the mean curvature of $F^k_0$ at $0$ is equal
to $\delta_k$. Hence the mean curvature of $F'_0$ vanishes at the
origin, which implies $F'_0$ is a single hyperplane or a pair of
parallel planes (by the constant rank of the second fundamental
form). The former case is ruled out by our definition of $r_{x,t}$
in (3.2), as it implies that the unit sphere is the largest tangent
sphere of $F'_0$ at the origin. The latter case was ruled out in
Corollary 6.1. Hence there is a $\delta>0$ such that $r_{x, t}H(x,
t)\ge \delta$ for all point $(x, t)\in\F$. \hfill$\square$

\vskip5pt

\noo{\bf Proof of Theorem 1.3}. If Theorem 1.3 is not true, there is
a sequence of normalized limit flows $\F'^k$  such that $r_k\to 0$,
where as in (3.2),
$$ r_k=
 \sup\{\rho:\ \B_\rho (z)\subset U_0'^k\ \text{and}\  |z-x|=\rho\}$$
is the radius of the largest ball which is contained in $U_0'^k$ and
tangential to $F_0'^k$ at $0$.

Regard $\F'^k$ as a hypersurface in the spacetime $\R^{n+1}\times
\R^1$, then it is a graph of a function $u_k$,
$$x_{n+2}=u_k(x_1, \cdots, x_{n+1}),$$
where $x_{n+2}=-t$, and $u_k$ satisfies the equation
$${\sum}_{i,j=1}^{n+1} (\delta_{ij}-\frac {u_iu_j}{|Du|^2})u_{ij} =1.$$
Since $\F'^k$ is normalized, the mean curvature of
$F'^k_0=\{u_k=0\}$ is equal to 1 at the origin, namely
$$|Du_k(0)|=1. $$
By Proposition 4.1 [20], $u_k$ is a convex function, and by Theorem
2.1 [20] (the version arXiv:math.DG/0404326), we have the estimate
$$u_k(x)\le C(1+|x|^2)$$
for some constant $C$ depending only on $n$. Hence $u_k$ converges
along a subsequence to a convex function $u$, whose level set is
still a solution to the mean curvature flow.

Since $r_k\to 0$ and $\F'^k$ is convex, the volume of the convex
sets $\{u_k<0\}\cap \B_R(0)$ converges to zero for any given $R>0$,
for otherwise $\F'^k$ is MCF of graph with uniformly bounded
gradient near the origin. Hence we have
$$u(0)=0,\ \ \ \ u\ge 0.$$
Recall that $|Du_k(0)|=1$ and $u_k$ is convex, by a rotation of
coordinates we assume that $\{x_{n+2}=x_1\}$ is a supporting plane
of $u_k$ at $0$, namely $u_k(x)\ge x_1$ for any $x\in\R^{n+1}$.
Hence we also have
$$u(x)\ge x_1\ \ \forall\ x\in\R^{n+1}. $$
By Lemma 2.6 of [20], the set $\{u=0\}$ is either a single point, or
a linear subspace of $\R^{n+1}$.

In the latter case, we may choose the axes such that $\{u=0\}$ is
the subspace spanned by the $x_{k+1}, \cdots, x_{n+1}$-axes. Then by
convexity, $u$ is independent of the variable $x_{k+1}, \cdots,
x_{n+1}$, namely $u(x)=u(x_1, \cdots, x_k)$. By restricting $u$ to
the subspace $\R^k$ spanned by the $x_1, \cdots, x_k$-axes, it
reduces to the former case, namely the set $\{u=0\}$ is a single
point.

Since the level set $F_t:=\{u=-t\}$ is a convex solution to the mean
curvature flow, by Huisken [8], it converges to a round point. Hence
we have $u(x)=O(|x|^2)$, which is in contradiction with the estimate
$u(x)\ge x_1$. Hence $r_k$ has a positive lower bound. This
completes the proof. \hfill $\square$

 \vskip5pt

\noo{\bf Proof of Corollary 1.3}. Assume to the contrary that (1.2)
is not true, then there exists $\delta_k\to 0$, and two sequences of
points $p_k$ and $q_k$, where $q_k\in Q_R(p_k)$, such that
$H(p_k)=\delta_k^{-1}$ and $H(q_k)\ge \delta_k^{-2}$.  Let
$\F^k=\D_{\delta_k^{-1}, p_k}(\F)$ be the normalized blow-up
sequence at $p_k$. Then $\F^k\cap \B_{2R}(0)$ converges smoothly to
a convex hypersurface. Hence $H_{\F^k}(\hat q_k)$ is uniformly
bounded, where $H_{\F^k}(\hat q_k)$ is the mean curvature of $\F^k$
at $\hat q_k=\D_{\delta_k^{-1}, p_k}(q_k)$.  But by assumption we
have $H_{\F^k}(\hat q_k)\ge \delta_k^{-1}$. We reach a
contradiction, as $\F^k$ converges locally smoothly.

Similarly one can prove the second inequality of (1.2).
\hfill$\square$

\vskip0pt

\section{Tangent flow}

First we state the following result, which is due to Huisken [10],
see also Ilmanen [15], and White [22].

\vskip5pt

\noo{\bf Proposition 8.1.} {\it A tangent flow to a mean convex flow
must be a shrinking cylinder, namely $\F'=S^m_r\times\R^{n-m}$ for
$1\le m\le n$; or a hyperplane $\F'=\R^n$,  where $S^m_r$ ($m\ge 1$)
denotes the shrinking sphere to the MCF in $\R^{m+1}$.}

\vskip5pt

From Theorems 1.1 and 1.2, we also have the follow result, which
means that the blow-up solution at a fixed first time singular point
is a shrinking sphere or cylinder.

\vskip5pt

\noo{\bf Corollary 8.1.} {\it A tangent flow at a (first time)
singular point cannot be a hyperplane.}

\noo{\bf Proof.}\ Let $(x_0, t_0)$ be a first time singular point.
For any $\tau>0$ small, let $x_\tau\in F_{t_0-\tau}$ be the point
closest to $x_0$. It suffices to prove that
 \begin{equation}
 c_1\tau^{1/2}\le r_\tau\le c_2\tau^{1/2}
 \end{equation}
where $r_\tau=\text{dist}(x_0, F_{t_0-\tau})$.

For the first inequality, consider the normalized blow-up sequence
$\F^\tau$ at $(x_\tau, t_-\tau)$, which converges along a
subsequence locally smoothly to a convex limit flow $\F'$. By the
$\kappa$-noncollapsing and since $\F^\tau$ is normalized blow-up
sequence, there is a sphere of radius $\kappa$, tangent to $F'_0$ at
the origin, contained in $U'_0$. Hence the limit flow $\F'=\{F'_t\}$
does not develop singularity for $t\in (0, c_0)$
($c_0=\big(\frac{\kappa}{2n}\big)^2$), which is equivalent by the
scaling (3.1) to that  $H(x_\tau, t_0-\tau)\tau^{1/2}\ge c_0$ (here
$H$ is the mean curvature of $\F$). Hence the velocity at $x_\tau$
of the MCF is greater than $c_0\tau^{-1/2}$, which implies that
$\frac {d}{d\tau} r_\tau\ge c_0\tau^{-1/2}$. Hence the distance
$|x_\tau-x_0|\ge c_0\tau^{1/2}$. We obtain the first inequality.

To prove the second inequality, we consider the sphere $\p
\B_{r_\tau}(x_0)$ at time $t_0-\tau$, at any given $\tau$. If
$r_\tau>\sqrt{2n\tau}$, the above sphere evolving by the MCF will
not shrink to the point $x_0$ at time $t_0$, and so $(x_0, t_0)$
cannot be a singular point of the MCF $\F$. $\hfill\square$


\noindent {\it Addresses}:

\noindent {\small Weimin Sheng: Department of Mathematics, Zhejiang
University, Hangzhou 310027,  China.

\vskip2pt

\noindent Xu-Jia Wang: Centre for Mathematics and Its Applications,
Australian National University, Canberra ACT 0200, Australia.

\vskip5pt

\noindent {\it E-mail}: weimins@zju.edu.cn, Xu-Jia.Wang@anu.edu.au}


\baselineskip=12.0pt
\parskip=0pt


\begin{thebibliography}{AB}

{\small
\bibitem [1]{B} K. Brakke,
          The motion of a surface by its mean curvature,
          Princeton Univ. Press, 1978.

\bibitem [2]{CGG} Y.G. Chen,  Y. Giga and S. Goto,
          Uniqueness and existence of viscosity solutions of generalized
          mean curvature flow equation, J. Diff. Geom. 33(1991), 749-786.



\bibitem [3]{Ec04} K. Ecker,
          Regularity theory for mean curvature flow,
          Birkhauser,  Boston,  2004.



\bibitem [4]{EH91} K. Ecker and G. Huisken,
          Interior estimates for hypersurfaces moving by mean curvature,
          Invent. Math. 105(1991), 547-569.



\bibitem [5]{ES91} L.C. Evans and J. Spruck,
          Motion of level sets by mean curvature,
          I, J. Diff. Geom. 33 (1991),  635-681;
          II, Trans. Amer. Math. Soc. 330 (1992), 321--332;
          III, J. Geom. Anal. 2 (1992), 121--150;
          IV, J. Geom. Anal. 5 (1995), 77--114.



\bibitem [6]{H95} R. Hamilton,
          The formation of singularities in the Ricci flow;
          Surveys in differential geometry, Vol. II,  7--136,
          International Press, 1995.



\bibitem [7]{H95-2} R. Hamilton,
          Harnack estimates for the mean curvature flow,
          J. Differential Geom., 41 (1995), 215-226.



\bibitem [8]{Hu84} G. Huisken,
          Flow by mean curvature of convex surfaces into spheres,
          J. Diff. Geom. 20 (1984),  237-266.



\bibitem [9]{Hu90} G. Huisken,
          Asymptotic behaviour for singularities of the mean curvature flow,
          J. Diff. Geom. 31 (1990), 285-299.



\bibitem [10]{Hu93} G. Huisken,
          Local and global behaviour of hypersurfaces moving by mean
          curvature, Proc. Symp. Pure Math. 54(1993), 175-191.



\bibitem [11]{HS99CV} G. Huisken and C. Sinestrari,
          Mean curvature flow singularities for mean convex surfaces,
          Calc. Var. PDE,  8 (1999), 1-14.



\bibitem [12]{HS99AC} G. Huisken and C. Sinestrari,
          Convexity estimates for mean curvature flow and singularities
          for mean convex surfaces, Acta Math. 183 (1999),  45-70.



\bibitem [13]{HS} G. Huisken and C. Sinestrari,
          Mean curvature flow with serguries,
          Inventiones Mathematicae, to appear.



\bibitem [14]{I1} T. Ilmanen,
          Elliptic regularization and partial regularity for motion by mean
          curvature,  Mem. Amer. Math. Soc.,  108 (1994), no. 520.



\bibitem [15]{I2} T. Ilmanen,
          Singularities of mean curvature flow of surfaces,
          preprint, 1995.



\bibitem [16]{L} G. Lieberman,
          Second order parabolic differential equations,
          World Scientific, Singapore, 1996.



\bibitem [17]{P1} G. Perelman,
          The entropy formula for the Ricci flow and its geometric applications,
          arXiv:math/0211159.



\bibitem [18]{P2} G. Perelman,
          Ricci flow with surgery on three manifolds,
          arXiv:math/0303109.



\bibitem [19]{Wa1} X.-J. Wang,
          Interior gradient estimates for mean curvature equations,
          Math. Z. 228 (1998), 73--81.



\bibitem [20]{Wa2} X.-J. Wang,
          Convex solutions to the mean curvature flow, \newline
          arXiv:math.DG/0404326.



\bibitem [21] {Wh2} B. White,
          The size of the singular set in mean curvature flow of mean-convex sets,
          J. Amer. Math. Soc. 13(2000), 665--695.



\bibitem [22] {Wh3} B. White,
          The nature of singularities in mean curvature flow of mean-convex sets,
          J. Amer. Math. Soc. 16(2003), 123--138.

\bibitem [23]{Z} X.-P. Zhu,
          Lectures on mean curvature flows.
          International Press, 2002.}

\end{thebibliography}
\enddocument